# Various applications of the (exponential) complete Bell polynomials

Donal F. Connon

16 January 2010

## Abstract

In a rather straightforward manner, we develop the well-known formula for the Stirling numbers of the first kind in terms of the (exponential) complete Bell polynomials where the arguments include the generalised harmonic numbers.

We also show how the (exponential) complete Bell polynomials feature in a number of other areas of mathematical interest.

## EXPLICIT FORMULA FOR THE STIRLING NUMBERS OF THE FIRST KIND

The Stirling numbers $s(n,k)$ of the first kind [23, p.56] are defined by the following generating function (the bracket symbol $\begin{bmatrix} n \\ k \end{bmatrix}$ is also employed)

$$(1.1) \qquad x(x-1)\cdots(x-n+1) = \sum_{k=0}^{n} s(n,k) x^k$$

which we may also express as the infinite series

$$(1.2) \qquad x(x-1)\cdots(x-n+1) = \sum_{k=0}^{\infty} s(n,k) x^k$$

where we define $s(n,k) = 0$ for $k \geq n+1$.

Letting $x \to -x$ in (1.2) results in

$$(1.3) \qquad x(x+1)\ldots(x+n-1) = \sum_{k=0}^{\infty} (-1)^{n+k} s(n,k) x^k = \sum_{k=1}^{\infty} (-1)^{n+k} s(n,k) x^k$$

We also consider the ascending factorial symbol $(x)_n$, also known as the Pochhamer symbol, defined by [23, p.16] as

$$(1.4) \qquad (x)_n = x(x+1)(x+2)\cdots(x+n-1) \text{ if } n > 0$$

$$(x)_0 = 1$$

The gamma function satisfies the well known recurrence relation

$$x\Gamma(x) = \Gamma(1+x)$$

and it is easily proved by induction that this may be extended to

(1.4.1) $$x(x+1)(x+2)\cdots(x+n-1)\Gamma(x) = \Gamma(n+x)$$

Hence we have

(1.5) $$(x)_n = \frac{\Gamma(n+x)}{\Gamma(x)} = \frac{x\Gamma(n+x)}{\Gamma(1+x)}$$

For the reasons discussed below, we shall find it more convenient to deal with the function $(x)_n$ divided by $x$. We consider the first derivative

(1.6) $$\frac{d}{dx}\frac{(x)_n}{x} = \frac{\Gamma(n+x)}{\Gamma(1+x)}[\psi(n+x) - \psi(1+x)]$$

where $\psi(x)$ is the digamma function defined by

(1.6.1) $$\psi(x) = \frac{d}{dx}\log\Gamma(x)$$

By definition we have

$$\psi(n+x) - \psi(1+x) = \frac{d}{dx}\log\Gamma(n+x) - \frac{d}{dx}\log\Gamma(1+x)$$

$$= \frac{d}{dx}\log\frac{\Gamma(n+x)}{\Gamma(1+x)}$$

$$= \frac{d}{dx}\log\frac{(x)_n}{x}$$

$$= \frac{d}{dx}\log(x)_n - \frac{1}{x}$$

Referring to the definition (1.4) of the Pochhamer symbol, we have

$$\log(x)_n = \log x + \log(x+1) + \cdots + \log(x+n-1)$$

and differentiation results in

$$\frac{(x)'_n}{(x)_n} = \frac{1}{x} + \frac{1}{x+1} + \cdots + \frac{1}{x+n-1} = H_n^{(1)}(x)$$

where $H_n^{(m)}(x)$ is the generalised harmonic number function defined by

(1.7) $$H_n^{(m)}(x) = \sum_{k=0}^{n-1} \frac{1}{(k+x)^m}$$

and we note that

$$H_n^{(m)}(1) = H_n^{(m)} = \sum_{k=1}^{n} \frac{1}{k^m}$$

Therefore we have the well known formula [23, p.14]

$$\psi(n+x) - \psi(1+x) = H_n^{(1)}(x) - \frac{1}{x}$$

$$= \sum_{k=1}^{n-1} \frac{1}{k+x}$$

and we denote $g(x) = \sum_{k=1}^{n-1} \frac{1}{k+x}$.

This gives us the relationship (which in turn leads us to consider the (exponential) complete Bell polynomials)

$$\frac{d}{dx}\frac{(x)_n}{x} = \frac{(x)_n}{x}\left[\psi(n+x) - \psi(1+x)\right]$$

(1.8) $$= \frac{(x)_n}{x} g(x)$$

We see that

$$g(0) = H_{n-1}^{(1)}$$

$$g^{(r)}(x) = (-1)^r r! \sum_{k=1}^{n-1} \frac{1}{(k+x)^{r+1}}$$

$$g^{(r)}(0) = (-1)^r r! H_{n-1}^{(r+1)}$$

As noted by Kölbig [14] and Coffey [4] we have

(1.9) $$\frac{d^r}{dx^r}e^{f(x)} = e^{f(x)}Y_r\left(f^{(1)}(x), f^{(2)}(x),..., f^{(r)}(x)\right)$$

where the (exponential) complete Bell polynomials may be defined by $Y_0 = 1$ and for $r \geq 1$

(1.10) $$Y_r(x_1,...,x_r) = \sum_{\pi(r)} \frac{r!}{k_1! k_2!... k_r!}\left(\frac{x_1}{1!}\right)^{k_1}\left(\frac{x_2}{2!}\right)^{k_2}\cdots\left(\frac{x_r}{r!}\right)^{k_r}$$

where the sum is taken over all partitions $\pi(r)$ of $r$, i.e. over all sets of integers $k_j$ such that

$$k_1 + 2k_2 + 3k_3 + \cdots + rk_r = r$$

The definition (1.10) immediately implies the following relation

(1.10.1) $$Y_r(ax_1, a^2x_2,..., a^r x_r) = a^r Y_r(x_1,...,x_r)$$

and with $a = 1$ we have

(1.10.2) $$Y_r(-x_1, x_2,..., (-1)^r x_r) = (-1)^r Y_r(x_1,...,x_r)$$

The complete Bell polynomials have integer coefficients and the first six are set out below (Comtet [6, p.307])

(1.11) $$Y_1(x_1) = x_1$$

$$Y_2(x_1, x_2) = x_1^2 + x_2$$

$$Y_3(x_1, x_2, x_3) = x_1^3 + 3x_1x_2 + x_3$$

$$Y_4(x_1, x_2, x_3, x_4) = x_1^4 + 6x_1^2x_2 + 4x_1x_3 + 3x_2^2 + x_4$$

$$Y_5(x_1, x_2, x_3, x_4, x_5) = x_1^5 + 10x_1^3x_2 + 10x_1^2x_3 + 15x_1x_2^2 + 5x_1x_4 + 10x_2x_3 + x_5$$

$$Y_6(x_1, x_2, x_3, x_4, x_5, x_6) = x_1^6 + 6x_1x_5 + 15x_2x_4 + 10x_3^2 + 15x_1^2x_4 + 15x_2^3 + 60x_1x_2x_3$$
$$+ 20x_1^3x_3 + 45x_1^2x_2^2 + 15x_1^4x_1 + x_6$$

The complete Bell polynomials are also given by the exponential generating function (Comtet [6, p.134])

$$(1.12) \qquad \exp\left(\sum_{j=1}^{\infty} x_j \frac{t^j}{j!}\right) = \sum_{n=0}^{\infty} Y_n(x_1,...,x_n) \frac{t^n}{n!}$$

Using (1.9) we see that

$$\frac{d^n}{dt^n} \exp\left(\sum_{j=1}^{\infty} x_j \frac{t^j}{j!}\right)\bigg|_{t=0} = Y_n(x_1,...,x_n)$$

and hence we note that (1.12) is simply the corresponding Maclaurin series.

We note that

$$\sum_{n=0}^{\infty} Y_n(ax_1,...,ax_n) \frac{t^n}{n!} = \exp\left(\sum_{j=1}^{\infty} ax_j \frac{t^j}{j!}\right) = \exp a\left(\sum_{j=1}^{\infty} x_j \frac{t^j}{j!}\right) = \left[\exp\left(\sum_{j=1}^{\infty} x_j \frac{t^j}{j!}\right)\right]^a$$

and thus we have

$$\left[\sum_{n=0}^{\infty} Y_n(x_1,...,x_n) \frac{t^n}{n!}\right]^a = \sum_{n=0}^{\infty} Y_n(ax_1,...,ax_n) \frac{t^n}{n!}$$

Let us now consider a function $f(t)$ which has a Taylor series expansion around $x$: we have

$$e^{f(x+t)} = \exp\left(\sum_{j=0}^{\infty} f^{(j)}(x) \frac{t^j}{j!}\right) = e^{f(x)} \exp\left(\sum_{j=1}^{\infty} f^{(j)}(x) \frac{t^j}{j!}\right)$$

$$= e^{f(x)} \left\{1 + \sum_{n=1}^{\infty} Y_n\left(f^{(1)}(x), f^{(2)}(x),..., f^{(n)}(x)\right) \frac{t^n}{n!}\right\}$$

We see that

$$\frac{d^m}{dx^m} e^{f(x)} = \frac{\partial^m}{\partial x^m} e^{f(x+t)}\bigg|_{t=0} = \frac{\partial^m}{\partial t^m} e^{f(x+t)}\bigg|_{t=0}$$

and we therefore obtain a derivation of (1.9) above

$$\frac{d^r}{dx^r} e^{f(x)} = e^{f(x)} Y_r\left(f^{(1)}(x), f^{(2)}(x),..., f^{(r)}(x)\right)$$

Suppose that $h'(x) = h(x)g(x)$ and let $f(x) = \log h(x)$. We see that

$$f'(x) = \frac{h'(x)}{h(x)} = g(x)$$

and then using (1.9) above we have

(1.13) $$\frac{d^r}{dx^r} h(x) = \frac{d^r}{dx^r} e^{\log h(x)} = h(x) Y_r\left(g(x), g^{(1)}(x), \ldots, g^{(r-1)}(x)\right)$$

In particular we have

$$\left.\frac{d^r}{dx^r} h(x)\right|_{x=0} = h(0) Y_r\left(g(0), g^{(1)}(0), \ldots, g^{(r-1)}(0)\right)$$

In our case, from (1.8) we have

$$h(x) = \frac{(x)_n}{x}$$

$$g(x) = \sum_{k=1}^{n-1} \frac{1}{k+x}$$

and it is easily seen that

(1.15) $$\left.\frac{d^r}{dx^r} \frac{(x)_n}{x}\right|_{x=0} = (n-1)! Y_r\left(H_{n-1}^{(1)}, -1! H_{n-1}^{(2)}, \ldots, (-1)^{r-1}(r-1)! H_{n-1}^{(r)}\right)$$

Differentiation of (1.3) also results in

(1.16) $$\frac{d^r}{dx^r} \frac{(x)_n}{x} = \sum_{k=0}^{\infty} (-1)^{n+k} (k-1) \cdots (k-r) s(n,k) x^{k-1-r}$$

and in particular we have

(1.17) $$\left.\frac{d^r}{dx^r} \frac{(x)_n}{x}\right|_{x=0} = (-1)^{n+r+1} r! s(n, r+1)$$

Therefore, equating (1.15) and (1.17), we obtain the known relationship for Stirling numbers of the first kind for $r \geq 0$

(1.18) $$s(n, r+1) = (-1)^{n+r+1} \frac{(n-1)!}{r!} Y_r\left(H_{n-1}^{(1)}, -1! H_{n-1}^{(2)}, \ldots, (-1)^{r-1}(r-1)! H_{n-1}^{(r)}\right)$$

The above relationship was previously derived by Kölbig [14] but this appears to be a more direct proof of this important formula.

The first few Stirling numbers $s(n,k)$ of the first kind are easily determined from (1.18); these are also reported in [21] and in the book by Srivastava and Choi [23, p.57]

(1.19)　　　$s(n,0) = \delta_{n,0}$

$$s(n,1) = (-1)^{n+1}(n-1)!$$

$$s(n,2) = (-1)^n (n-1)! H^{(1)}_{n-1}$$

$$s(n,3) = (-1)^{n+1} \frac{(n-1)!}{2}\left\{\left(H^{(1)}_{n-1}\right)^2 - H^{(2)}_{n-1}\right\}$$

$$s(n,4) = (-1)^n \frac{(n-1)!}{6}\left\{\left(H^{(1)}_{n-1}\right)^3 - 3H^{(1)}_{n-1}H^{(2)}_{n-1} + 2H^{(3)}_{n-1}\right\}$$

□

Using (1.18) gives us

$$(-1)^{n+r-k+1}\frac{(r-k)!}{(n-1)!}s(n,r-k+1) = Y_{r-k}\left(H^{(1)}_{n-1}, -1!H^{(2)}_{n-1}, \ldots, (-1)^{r-k-1}(r-k-1)!H^{(r-k)}_{n-1}\right)$$

and since [5, p.415]

(1.20)　　　$Y_{r+1}(x_1,\ldots,x_{r+1}) = \sum_{k=0}^{r}\binom{r}{k}Y_{r-k}(x_1,\ldots,x_{r-k})x_{k+1}$

we obtain the recurrence formula given by Shen [21]

$$(r+1)s(n,r+2) = -\sum_{k=0}^{r} s(n,r-k+1)H^{(k+1)}_{n-1}$$

□

Let us now see why we initially directed our attention to $(x)_n$ divided by $x$. Differentiation of (1.5) results in

(1.21)　　　$\dfrac{d^r}{dx^r}(x)_n = \dfrac{d^r}{dx^r}\dfrac{\Gamma(n+x)}{\Gamma(x)} = \sum_{k=0}^{\infty}(-1)^{n+k}k(k-1)\cdots(k-r+1)s(n,k)x^{k-r}$

and in particular we have

(1.22) $$\left.\frac{d^r}{dx^r}(x)_n\right|_{x=0} = (-1)^{n+r} r! s(n,r)$$

Alternatively we have the first derivative

$$\frac{d}{dx}(x)_n = \frac{d}{dx}\frac{\Gamma(n+x)}{\Gamma(x)} = \frac{\Gamma(x)\Gamma'(n+x) - \Gamma'(x)\Gamma(n+x)}{\Gamma^2(x)}$$

$$= \frac{\Gamma(n+x)}{\Gamma(x)}[\psi(n+x) - \psi(x)]$$

Therefore we have

(1.23) $$\frac{d}{dx}(x)_n = (x)_n H_n^{(1)}(x) = \frac{(x)_n}{x} x H_n^{(1)}(x)$$

and it is easily seen from (1.7) that

$$\lim_{x \to 0} x H_n^{(1)}(x) = 1$$

Hence we obtain

$$\left.\frac{d}{dx}(x)_n\right|_{x=0} = (n-1)!$$

and using (1.19) we see again that

$$s(n,1) = (-1)^{n+1}(n-1)!$$

It should be noted that $\left.\frac{d^r}{dx^r}(x)_n\right|_{x=0}$ may also be determined directly; however, the mathematics is somewhat more cumbersome. For example, let us consider $\left.\frac{d^2}{dx^2}(x)_n\right|_{x=0}$.

Differentiating (1.23) gives us

$$\frac{d^2}{dx^2}(x)_n = (x)_n\left(\left[H_n^{(1)}(x)\right]^2 - H_n^{(2)}(x)\right)$$

We have by definition

$$H_n^{(1)}(x+1) = \frac{1}{x+1} + \frac{1}{x+2} + \cdots + \frac{1}{x+n}$$

and therefore

$$H_n^{(1)}(x) = \frac{1}{x} + \frac{1}{x+1} + \cdots + \frac{1}{x+n-1} = \frac{1}{x} - \frac{1}{x+n} + H_n^{(1)}(x+1)$$

$$= \frac{n}{x(x+n)} + H_n^{(1)}(x+1)$$

This gives us

$$\left[H_n^{(1)}(x)\right]^2 = \frac{n^2}{x^2(x+n)^2} + \left[H_n^{(1)}(x+1)\right]^2 + \frac{2n}{x(x+n)} H_n^{(1)}(x+1)$$

Similarly we have

$$H_n^{(2)}(x) = \frac{1}{x^2} - \frac{1}{(x+n)^2} + H_n^{(2)}(x+1)$$

$$= \frac{2nx + n^2}{x^2(x+n)^2} + H_n^{(2)}(x+1)$$

and we easily see that

$$\left[H_n^{(1)}(x)\right]^2 - H_n^{(2)}(x) = \left[H_n^{(1)}(x+1)\right]^2 - H_n^{(2)}(x+1) + \frac{2n}{x(x+n)} H_n^{(1)}(x+1) - \frac{2nx}{x^2(x+n)^2}$$

This gives us

$$\lim_{x \to 0} x \left( \left[H_n^{(1)}(x)\right]^2 - H_n^{(2)}(x) \right) = 2 H_{n-1}^{(1)}$$

and we then obtain

$$\frac{d^2}{dx^2}(x)_n \bigg|_{x=0} = \lim_{x \to 0} (x)_n \left( \left[H_n^{(1)}(x)\right]^2 - H_n^{(2)}(x) \right)$$

$$= \lim_{x \to 0} \frac{(x)_n}{x} x \left( \left[H_n^{(1)}(x)\right]^2 - H_n^{(2)}(x) \right)$$

$$= 2(n-1)!H_{n-1}^{(1)}$$

Hence we have using (1.22)

$$s(n,2) = (-1)^n (n-1)!H_{n-1}^{(1)}$$

□

From (1.23) we find that

$$\frac{d^r}{dx^r}(x)_n = (x)_n Y_r\left(H_n^{(1)}(x), -1!H_n^{(2)}(x), \ldots, (-1)^{r-1}(r-1)!H_n^{(r)}(x)\right)$$

and using (1.21) we see that

$$(-1)^{n+r} r! s(n,r) = \lim_{x \to 0} (x)_n Y_r\left(H_n^{(1)}(x), -1!H_n^{(2)}(x), \ldots, (-1)^{r-1}(r-1)!H_n^{(r)}(x)\right)$$

## CAUCHY'S GENERATING FUNCTION FOR THE STIRLING NUMBERS OF THE FIRST KIND

Using the binomial theorem we have for $|z| < 1$

(2.1) $$\frac{1}{(1-z)^x} = \sum_{n=0}^{\infty} \frac{(x)_n}{n!} z^n$$

and, using (1.22), differentiation with respect to $x$ results in

(2.2) $$-\frac{\log(1-z)}{(1-z)^x} = \sum_{n=0}^{\infty} \frac{(x)_n H_n^{(1)}(x)}{n!} z^n$$

More generally we have

$$(-1)^r \frac{\log^r(1-z)}{(1-z)^x} = \sum_{n=0}^{\infty} \frac{d^r}{dx^r}(x)_n \frac{z^n}{n!}$$

and letting $x = 0$ we see that

$$(-1)^r \log^r(1-z) = \sum_{n=0}^{\infty} \frac{d^r}{dx^r}(x)_n \bigg|_{x=0} \frac{z^n}{n!}$$

Then using (1.22)

$$\left.\frac{d^r}{dx^r}(x)_n\right|_{x=0} = (-1)^{n+r} r! s(n,r)$$

we then easily obtain the well-known Maclaurin expansion due to Cauchy [23, p.56]

(2.3) $$\log^r(1-z) = r! \sum_{n=0}^{\infty} (-1)^n s(n,r) \frac{z^n}{n!}$$

Since $s(n,r) = 0$ for $r \geq n+1$ this may be expressed as

$$\log^r(1-z) = r! \sum_{n=r}^{\infty} (-1)^n s(n,r) \frac{z^n}{n!}$$

A different proof of (2.3) was given by Póyla and Szegö in [18, p.227].

Letting $x = 1$ in (2.2) we obtain the well-known generating function for the harmonic numbers

(2.3.1) $$-\frac{\log(1-z)}{1-z} = \sum_{n=0}^{\infty} H_n^{(1)} z^n$$

Differentiating (2.2) gives us

(2.4) $$\frac{\log^2(1-z)}{(1-z)^x} = \sum_{n=0}^{\infty} \frac{(x)_n \left( \left[H_n^{(1)}(x)\right]^2 - H_n^{(2)}(x) \right)}{n!} z^n$$

and with $x = 1$ we have

(2.5) $$\frac{\log^2(1-z)}{1-z} = \sum_{n=0}^{\infty} \left( \left[H_n^{(1)}\right]^2 - H_n^{(2)} \right) z^n$$

Such series are considered in more detail in, for example, [7] and the references contained therein.

Integrating (2.5) results in

$$-\frac{1}{3}\log^3(1-z) = \sum_{n=0}^{\infty} \left( \left[H_n^{(1)}\right]^2 - H_n^{(2)} \right) \frac{z^{n+1}}{n+1}$$

$$= \sum_{n=0}^{\infty} \left( \left[H_{n+1}^{(1)} - \frac{1}{n+1}\right]^2 - H_{n+1}^{(2)} + \frac{1}{(n+1)^2} \right) \frac{z^{n+1}}{n+1}$$

$$= \sum_{n=0}^{\infty}\left(\left[H_{n+1}^{(1)}\right]^2 - \frac{2H_{n+1}^{(1)}}{n+1} - H_{n+1}^{(2)} + \frac{2}{(n+1)^2}\right)\frac{z^{n+1}}{n+1}$$

$$= \sum_{n=1}^{\infty}\left(\left[H_{n}^{(1)}\right]^2 - \frac{2H_{n}^{(1)}}{n} - H_{n}^{(2)} + \frac{2}{n^2}\right)\frac{z^{n}}{n}$$

We therefore have

(2.6) $$-\frac{1}{3}\log^3(1-z) = \sum_{n=1}^{\infty}\left(\left[H_n^{(1)}\right]^2 - H_n^{(2)}\right)\frac{z^n}{n} - 2\sum_{n=1}^{\infty}\frac{H_n^{(1)}}{n^2}z^n + 2Li_3(z)$$

where $Li_3(z)$ is the polylogarithm defined by

(2.7) $$Li_s(z) = \sum_{n=1}^{\infty}\frac{z^n}{n^s}$$

From (2.3) we have with $r = 3$

$$\log^3(1-z) = 6\sum_{n=0}^{\infty}(-1)^n s(n,3)\frac{z^n}{n!}$$

$$= -3\sum_{n=1}^{\infty}\left(\left[H_{n-1}^{(1)}\right]^2 - H_{n-1}^{(2)}\right)\frac{z^n}{n}$$

$$= -3\sum_{n=1}^{\infty}\left(\left[H_n^{(1)}\right]^2 - \frac{2H_n^{(1)}}{n} - H_n^{(2)} + \frac{2}{n^2}\right)\frac{z^n}{n}$$

in agreement with (2.6) above.

## COPPO'S FORMULA

We have the well-known partial fraction decomposition [13, p.188]

(3.1) $$f(x) = \frac{n!}{x(x+1)\ldots(x+n)} = \sum_{k=0}^{n}\binom{n}{k}\frac{(-1)^k}{k+x}$$

$$= x^{-1}\binom{n+x}{n}^{-1}, \quad x \notin \{0,-1,\ldots,-n\}$$

and we note from (1.4.1) that

$$\frac{n!}{x(1+x)...(n+x)} = \frac{\Gamma(n+1)\Gamma(x)}{\Gamma(n+1+x)}$$

or equivalently

$$\frac{n!}{x(1+x)...(n+x)} = \frac{\Gamma(n+1)}{(x)_{n+1}}$$

Differentiation results in

$$f'(x) = -\frac{\Gamma(n+1)}{(x)_{n+1}} H_{n+1}^{(1)}(x)$$

and referring to (1.13) we therefore have the $r$ th derivative

$$f^{(r)}(x) = \frac{\Gamma(n+1)}{(x)_{n+1}} Y_r\left(-H_{n+1}^{(1)}(x), 1!H_{n+1}^{(2)}(x), ..., (-1)^r(r-1)!H_{n+1}^{(r)}(x)\right)$$

We also have from (3.1)

$$f^{(r)}(x) = (-1)^r r! \sum_{k=0}^{n} \binom{n}{k} \frac{(-1)^k}{(k+x)^{r+1}}$$

and we therefore obtain

$$(3.2) \quad \sum_{k=0}^{n} \binom{n}{k} \frac{(-1)^k}{(k+x)^{r+1}} = \frac{\Gamma(n+1)}{(x)_{n+1}} \frac{(-1)^r}{r!} Y_r\left(-0!H_{n+1}^{(1)}(x), 1!H_{n+1}^{(2)}(x), ..., (-1)^r(r-1)!H_{n+1}^{(r)}(x)\right)$$

It may be noted that Coppo [10] has expressed this in a slightly different form

$$(3.3) \quad \sum_{k=0}^{n} \binom{n}{k} \frac{(-1)^k}{(k+x)^{r+1}} = \frac{\Gamma(n+1)}{(x)_{n+1}} \frac{1}{r!} Y_r\left(0!H_{n+1}^{(1)}(x), 1!H_{n+1}^{(2)}(x), ..., (r-1)!H_{n+1}^{(r)}(x)\right)$$

and reference to (1.10.2) shows that these are equivalent statements. It may be noted that reference to the right-hand side of (3.3) shows that $\sum_{k=0}^{n} \binom{n}{k} \frac{(-1)^k}{(k+x)^{r+1}}$ is positive for $x > 0$ (which otherwise does not appear to be immediately obvious).

Particular cases of Coppo's formula are set out below.

(3.4) $$\sum_{k=0}^{n}\binom{n}{k}\frac{(-1)^k}{(k+x)} = \frac{\Gamma(n+1)\Gamma(x)}{\Gamma(n+1+x)}Y_0 = \frac{\Gamma(n+1)\Gamma(x)}{\Gamma(n+1+x)}$$

(3.5) $$\sum_{k=0}^{n}\binom{n}{k}\frac{(-1)^k}{(k+x)^2} = -\frac{\Gamma(n+1)\Gamma(x)}{\Gamma(n+1+x)}Y_1\left(-H_{n+1}^{(1)}(x)\right) = \frac{\Gamma(n+1)\Gamma(x)}{\Gamma(n+1+x)}H_{n+1}^{(1)}(x)$$

$$\sum_{k=0}^{n}\binom{n}{k}\frac{(-1)^k}{(k+x)^3} = \frac{\Gamma(n+1)\Gamma(x)}{\Gamma(n+1+x)}\frac{1}{2}Y_2\left(-H_{n+1}^{(1)}(x), H_{n+1}^{(2)}(x)\right)$$

(3.6) $$= \frac{\Gamma(n+1)\Gamma(x)}{\Gamma(n+1+x)}\frac{1}{2}\left(\left[H_{n+1}^{(1)}(x)\right]^2 + H_{n+1}^{(2)}(x)\right)$$

$$\sum_{k=0}^{n}\binom{n}{k}\frac{(-1)^k}{(k+x)^4} = -\frac{\Gamma(n+1)\Gamma(x)}{\Gamma(n+1+x)}\frac{1}{6}Y_3\left(-H_{n+1}^{(1)}(x), H_{n+1}^{(2)}(x), -2H_{n+1}^{(3)}(x)\right)$$

(3.7) $$= \frac{\Gamma(n+1)\Gamma(x)}{\Gamma(n+1+x)}\frac{1}{6}\left(\left[H_{n+1}^{(1)}(x)\right]^3 + 3H_{n+1}^{(1)}(x)H_{n+1}^{(2)}(x) + 2H_{n+1}^{(3)}(x)\right)$$

$\square$

We now wish to consider representations for $\sum_{k=1}^{n}\binom{n}{k}\frac{(-1)^k}{k^r}$. We see that

$$\sum_{k=1}^{n}\binom{n}{k}\frac{(-1)^k}{k+x} = \sum_{k=0}^{n}\binom{n}{k}\frac{(-1)^k}{k+x} - \frac{1}{x} = \frac{\Gamma(n+1)\Gamma(x)}{\Gamma(n+1+x)} - \frac{1}{x}$$

$$= \frac{\Gamma(n+1)x\Gamma(x) - \Gamma(n+1+x)}{x\Gamma(n+1+x)}$$

$$= \frac{\Gamma(n+1)\Gamma(1+x) - \Gamma(n+1+x)}{x\Gamma(n+1+x)}$$

and we have the limit as $x \to 0$

$$\sum_{k=1}^{n}\binom{n}{k}\frac{(-1)^k}{k} = \lim_{x \to 0}\frac{\Gamma(n+1)x\Gamma(x) - \Gamma(n+1+x)}{x\Gamma(n+1+x)} \approx \frac{0}{0}$$

Applying L'Hôpital's rule we obtain

$$= \lim_{x \to 0} \frac{\Gamma(n+1)\Gamma'(1+x) - \Gamma'(n+1+x)}{x\Gamma'(n+1+x) + \Gamma(n+1+x)}$$

$$= \frac{\Gamma(n+1)\Gamma'(1) - \Gamma'(n+1)}{\Gamma(n+1)}$$

$$= \psi(1) - \psi(n+1)$$

$$= -H_n^{(1)}$$

Hence we obtain Euler's well-known identity

(3.8) $\quad -\sum_{k=1}^{n} \binom{n}{k} \frac{(-1)^k}{k} = H_n^{(1)}$

Similarly we also have

$$\sum_{k=1}^{n} \binom{n}{k} \frac{(-1)^k}{(k+x)^2} = \frac{\Gamma(n+1)\Gamma(x)}{\Gamma(n+1+x)} H_{n+1}^{(1)}(x) - \frac{1}{x^2}$$

$$= \frac{\Gamma(n+1)x^2\Gamma(x)H_{n+1}^{(1)}(x) - \Gamma(n+1+x)}{x^2\Gamma(n+1+x)}$$

$$= \frac{\Gamma(n+1)\Gamma(1+x)x H_{n+1}^{(1)}(x) - \Gamma(n+1+x)}{x^2\Gamma(n+1+x)}$$

and since $\lim_{x \to 0} x H_n^{(1)}(x) = 1$ we see that we may apply L'Hôpital's rule to obtain

$$\sum_{k=1}^{n} \binom{n}{k} \frac{(-1)^k}{k^2}$$

$$= \lim_{x \to 0} \frac{\Gamma(n+1)\Gamma(1+x)[H_{n+1}^{(1)}(x) - x H_{n+1}^{(2)}(x)] + \Gamma(n+1)\Gamma'(1+x)x H_{n+1}^{(1)}(x) - \Gamma'(n+1+x)}{x^2\Gamma'(n+1+x) + 2x\Gamma(n+1+x)}$$

Since $\lim_{x \to 0}[H_{n+1}^{(1)}(x) - x H_{n+1}^{(2)}(x)] = H_n^{(1)}$ we have

$$\approx \frac{\Gamma(n+1)H_n^{(1)} + \Gamma(n+1)\Gamma'(1) - \Gamma'(n+1)}{0}$$

$$\approx \frac{\Gamma(n+1)\left[H_n^{(1)} + \psi(1) - \psi(n+1)\right]}{0}$$

$$\approx \frac{0}{0}$$

Hence we may apply L'Hôpital's rule again to obtain

$$= \lim_{x \to 0} \frac{\Gamma(n+1)\Gamma(1+x)[-2H_{n+1}^{(2)}(x) + 2x\,H_{n+1}^{(3)}(x)]}{x^2\Gamma''(n+1+x) + 2x\Gamma'(n+1+x) + 2x\Gamma'(n+1+x) + 2\Gamma(n+1+x)}$$

$$+ \lim_{x \to 0} \frac{2\Gamma(n+1)\Gamma'(1+x)[H_{n+1}^{(1)}(x) - x\,H_{n+1}^{(2)}(x)]}{x^2\Gamma''(n+1+x) + 2x\Gamma'(n+1+x) + 2x\Gamma'(n+1+x) + 2\Gamma(n+1+x)}$$

$$+ \lim_{x \to 0} \frac{\Gamma(n+1)\Gamma''(1+x)x\,H_{n+1}^{(1)}(x)}{x^2\Gamma''(n+1+x) + 2x\Gamma'(n+1+x) + 2x\Gamma'(n+1+x) + 2\Gamma(n+1+x)}$$

$$- \lim_{x \to 0} \frac{\Gamma''(n+1+x)}{x^2\Gamma''(n+1+x) + 2x\Gamma'(n+1+x) + 2x\Gamma'(n+1+x) + 2\Gamma(n+1+x)}$$

$$= \frac{1}{2}[-2H_n^{(2)} + 2\Gamma'(1)H_n^{(1)}] - \frac{1}{2}\left[\Gamma''(1) - \frac{\Gamma''(n+1)}{\Gamma(n+1)}\right] = L$$

Differentiation of (1.6.1) gives us

$$\psi'(x) = \frac{\Gamma(x)\Gamma''(x) - (\Gamma'(x))^2}{\Gamma^2(x)} = \varsigma(2, x)$$

so that

$$\frac{\Gamma''(x)}{\Gamma(x)} = \psi^2(x) + \varsigma(2, x)$$

where $\varsigma(s, x)$ is the Hurwitz zeta function.

Therefore we have

$$\frac{\Gamma''(n+1)}{\Gamma(n+1)} = \psi^2(n+1) + \varsigma(2, n+1)$$

$$\Gamma''(1) = \psi^2(1) + \varsigma(2)$$

and so

$$\Gamma''(1) - \frac{\Gamma''(n+1)}{\Gamma(n+1)} = \varsigma(2) - \varsigma(2, n+1) + \psi^2(1) - \psi^2(n+1)$$

From the definition of the Hurwitz zeta function it readily follows that

$$\varsigma(s, n+x) = \varsigma(s, x) - H_n^{(s)}(x)$$

and in particular we have

$$\varsigma(2) - \varsigma(2, n+1) = H_n^{(2)}$$

Since $\psi(n+1) - \psi(1) = H_n^{(1)}$ we have

$$\psi^2(n+1) + \psi^2(1) - 2\psi(1)\psi(n+1) = \left[H_n^{(1)}\right]^2$$

Therefore we obtain

$$\Gamma''(1) - \frac{\Gamma''(n+1)}{\Gamma(n+1)} = H_n^{(2)} - \left[H_n^{(1)}\right]^2 + 2\psi^2(1) - 2\psi(1)\psi(n+1)$$

and thus

$$L = \frac{1}{2}\left(-H_n^{(2)} - \left[H_n^{(1)}\right]^2 + 2\psi^2(1) - 2\psi(1)\psi(n+1) + 2\psi(1)H_n^{(1)}\right)$$

which simplifies to

$$L = -\frac{1}{2}\left(H_n^{(2)} + \left[H_n^{(1)}\right]^2\right)$$

so that we obtain

(3.9) $$\sum_{k=1}^{n}\binom{n}{k}\frac{(-1)^k}{k^2} = -\frac{1}{2}\left(H_n^{(2)} + \left[H_n^{(1)}\right]^2\right)$$

The above result is well known; see for example the paper by Flajolet and Sedgewick [12] where they also reported the following identity

$$-S_n(3) = -\sum_{k=1}^{n}\binom{n}{k}\frac{(-1)^k}{k^3} = \frac{1}{6}\left(H_n^{(1)}\right)^3 + \frac{1}{2}H_n^{(1)}H_n^{(2)} + \frac{1}{3}H_n^{(3)}$$

which may also be derived by the above limiting process (albeit with much more tedious labour).

Defining $S_n(r)$ by

(3.10) $$S_n(r) = \sum_{k=1}^{n} \binom{n}{k} \frac{(-1)^k}{k^r}$$

then Flajolet and Sedgewick [12] showed that $S_n(r)$ can be expressed in terms of the generalised harmonic numbers as

(3.11) $$-S_n(r) = \sum_{1m_1+2m_2+3m_3\ldots=r} \frac{1}{m_1!\, m_2!\, m_3!\ldots m_r!} \left(\frac{H_n^{(1)}}{1}\right)^{m_1} \left(\frac{H_n^{(2)}}{2}\right)^{m_2} \left(\frac{H_n^{(3)}}{3}\right)^{m_3} \ldots \left(\frac{H_n^{(r)}}{r}\right)^{m_r}$$

Referring to (1.10) we then note that this may be written as

$$-S_n(r) = \frac{1}{r!} Y_r\left(0!H_n^{(1)}, 1!H_n^{(2)}, \ldots, (r-1)!H_n^{(r)}\right)$$

and hence we have

(3.12) $$S_n(r) = \sum_{k=1}^{n} \binom{n}{k} \frac{(-1)^k}{k^r} = -\frac{1}{r!} Y_r\left(0!H_n^{(1)}, 1!H_n^{(2)}, \ldots, (r-1)!H_n^{(r)}\right)$$

which may be contrasted with (3.2).

Rather belatedly, I noted that letting $x \to 1+x$ in (3.2) gives us

$$\sum_{k=0}^{n} \binom{n}{k} \frac{(-1)^k}{(k+1+x)^{r+1}} = \frac{\Gamma(n+1)}{(1+x)_{n+1}} \frac{(-1)^r}{r!} Y_r\left(-0!H_{n+1}^{(1)}(1+x), 1!H_{n+1}^{(2)}(1+x), \ldots, (-1)^r(r-1)!H_{n+1}^{(r)}(1+x)\right)$$

whereupon letting $x = 0$ results in

$$\sum_{k=0}^{n} \binom{n}{k} \frac{(-1)^k}{(k+1)^{r+1}} = \frac{1}{n+1} \frac{(-1)^r}{r!} Y_r\left(-0!H_{n+1}^{(1)}, 1!H_{n+1}^{(2)}, \ldots, (-1)^r(r-1)!H_{n+1}^{(r)}\right)$$

With $n \to n-1$ this becomes

$$\sum_{k=0}^{n-1} \binom{n-1}{k} \frac{(-1)^k}{(k+1)^{r+1}} = \frac{1}{n} \frac{(-1)^r}{r!} Y_r\left(-0!H_n^{(1)}, 1!H_n^{(2)}, \ldots, (-1)^r(r-1)!H_n^{(r)}\right)$$

Reindexing gives us

$$\sum_{m=1}^{n}\binom{n-1}{m-1}\frac{(-1)^{m+1}}{m^{r+1}} = \frac{1}{n}\frac{(-1)^r}{r!}Y_r\left(-0!H_n^{(1)}, 1!H_n^{(2)}, \ldots, (-1)^r(r-1)!H_n^{(r)}\right)$$

Since $\dfrac{1}{m}\binom{n-1}{m-1} = \dfrac{1}{n}\binom{n}{m}$ this becomes

(3.13) $$\sum_{m=1}^{n}\binom{n}{m}\frac{(-1)^{m+1}}{m^{r}} = \frac{(-1)^r}{r!}Y_r\left(-0!H_n^{(1)}, 1!H_n^{(2)}, \ldots, (-1)^r(r-1)!H_n^{(r)}\right)$$

and having regard to (1.10.2) we see that

(3.14) $$Y_r\left(0!H_n^{(1)}, 1!H_n^{(2)}, \ldots, (r-1)!H_n^{(r)}\right) = (-1)^r Y_r\left(-0!H_n^{(1)}, 1!H_n^{(2)}, \ldots, (-1)^r(r-1)!H_n^{(r)}\right)$$

Hence we obtain as before

(3.15) $$\sum_{m=1}^{n}\binom{n}{m}\frac{(-1)^{m+1}}{m^{r}} = \frac{1}{r!}Y_r\left(0!H_n^{(1)}, 1!H_n^{(2)}, \ldots, (r-1)!H_n^{(r)}\right)$$

□

We now multiply (3.12) by $x^r$ and make the summation to obtain

(3.16) $$\sum_{r=0}^{\infty} x^r \sum_{k=1}^{n}\binom{n}{k}\frac{(-1)^{k+1}}{k^{r}} = \sum_{r=0}^{\infty} Y_r\left(0!H_n^{(1)}, 1!H_n^{(2)}, \ldots, (r-1)!H_n^{(r)}\right)\frac{x^r}{r!}$$

and comparing this with (1.12) we have

(3.17) $$\sum_{r=0}^{\infty} x^r \sum_{k=1}^{n}\binom{n}{k}\frac{(-1)^{k+1}}{k^{r}} = \exp\left(\sum_{r=1}^{\infty}\frac{H_n^{(r)}}{r}x^r\right)$$

Let

$$f(x) = \log\Gamma(n+x) - \log\Gamma(1+x)$$

Then we have

$$f^{(m)}(0) = \psi^{(m-1)}(n) - \psi^{(m-1)}(1) = (-1)^{m+1}(m-1)!H_{n-1}^{(m)}$$

and we therefore have the Maclaurin expansion

(3.17) $$f(x) = \log\Gamma(n+x) - \log\Gamma(1+x) = \log\Gamma(n) + \sum_{m=1}^{\infty}\frac{(-1)^{m+1}}{m}H_{n-1}^{(m)}x^m$$

This may be written as

(3.18) $$\frac{\Gamma(n+x)}{\Gamma(1+x)\Gamma(n)} = \exp\left[\sum_{m=1}^{\infty} \frac{(-1)^{m+1}}{m} H_{n-1}^{(m)} x^m\right]$$

as previously noted by Wilf [24]. We also note that

$$\frac{\Gamma(1+x)\Gamma(n)}{\Gamma(n+x)} = \exp\left[\sum_{m=1}^{\infty} \frac{(-1)^m}{m} H_{n-1}^{(m)} x^m\right]$$

With $x \to -x$ and $n \to n+1$ (3.18) becomes

$$\frac{\Gamma(n+1-x)}{\Gamma(1-x)\Gamma(n+1)} = \exp\left[-\sum_{m=1}^{\infty} \frac{H_n^{(m)}}{m} x^m\right]$$

and therefore we have

$$\frac{\Gamma(1-x)\Gamma(n+1)}{\Gamma(n+1-x)} = \exp\left[\sum_{m=1}^{\infty} \frac{H_n^{(m)}}{m} x^m\right]$$

Hence we obtain

(3.19) $$\frac{\Gamma(1-x)\Gamma(n+1)}{\Gamma(n+1-x)} = \sum_{r=0}^{\infty} x^r \sum_{k=1}^{n} \binom{n}{k} \frac{(-1)^{k+1}}{k^r}$$

From (3.16) we have

$$\log \Gamma(n+x) - \log \Gamma(1+x) = \log \Gamma(n) + \sum_{m=1}^{\infty} \frac{(-1)^{m+1}}{m} H_{n-1}^{(m)} x^m$$

With $x \to -x$ and $n \to n+1$ this becomes

$$\log \frac{\Gamma(1-x)\Gamma(n+1)}{\Gamma(n+1-x)} = \sum_{m=1}^{\infty} \frac{H_n^{(m)}}{m} x^m$$

Therefore using (6.3) we have

$$\frac{\Gamma(1-x)\Gamma(n+1)}{\Gamma(n+1-x)} = \sum_{m=0}^{\infty} Y_m\left(H_n^{(1)}, 1! H_n^{(2)}, \ldots, (m-1)! H_n^{(m)}\right) \frac{x^m}{m!}$$

Employing (1.12) we see that

$$\sum_{m=0}^{\infty} Y_m\left(H_n^{(1)}, 1!H_n^{(2)}, \ldots, (m-1)!H_n^{(m)}\right)\frac{x^m}{m!} = \exp\left(\sum_{m=1}^{\infty}\frac{H_n^{(m)}}{m}x^m\right)$$

and hence we have

(3.20) $$\frac{\Gamma(1-x)\Gamma(n+1)}{\Gamma(n+1-x)} = \sum_{r=0}^{\infty} x^r \sum_{k=1}^{n}\binom{n}{k}\frac{(-1)^{k+1}}{k^r}$$

Using (1.12) and (3.2) gives us a more general identity

$$\sum_{r=0}^{\infty} t^r \frac{(-1)^r (x)_{n+1}}{\Gamma(n+1)} \sum_{k=0}^{n}\binom{n}{k}\frac{(-1)^k}{(k+x)^{r+1}} = \sum_{r=0}^{\infty} Y_r\left(-0!H_{n+1}^{(1)}(x), 1!H_{n+1}^{(2)}(x), \ldots, (-1)^r(r-1)!H_{n+1}^{(r)}(x)\right)\frac{t^r}{r!}$$

$$= \exp\left(\sum_{r=1}^{\infty}(-1)^r \frac{H_{n+1}^{(r)}(x)}{r}t^r\right)$$

In equation (3.67a) of [7] we showed that

(3.21) $$\sum_{n=1}^{\infty} t^n \sum_{k=1}^{n}\binom{n}{k}\frac{(-1)^{k+1}}{k^r} = \frac{1}{t-1}Li_r\left(\frac{-t}{1-t}\right)$$

and using (3.12) we see that

(3.22) $$\frac{1}{t-1}Li_r\left(\frac{-t}{1-t}\right) = \frac{1}{r!}\sum_{n=1}^{\infty} Y_r\left(0!H_n^{(1)}, 1!H_n^{(2)}, \ldots, (r-1)!H_n^{(r)}\right)t^n$$

where $Li_r(x)$ is the polylogarithm function defined by

$$Li_r(x) = \sum_{k=1}^{\infty}\frac{x^k}{k^r}$$

We note that

$$Li_r(-1) = \sum_{k=1}^{\infty}\frac{(-1)^k}{k^r} = -\varsigma_s(r)$$

where $\varsigma_s(r)$ is the alternating Riemann zeta function and

$$\varsigma_s(r) = (1-2^{1-r})\varsigma(r)$$

Therefore with $t = 1/2$ in (3.22) we obtain

$$-2Li_r(-1) = \frac{1}{r!}\sum_{n=1}^{\infty}\frac{1}{2^n}Y_r\left(0!H_n^{(1)}, 1!H_n^{(2)}, \ldots, (r-1)!H_n^{(r)}\right)$$

and thus we obtain

(3.23) $$\varsigma(r) = \frac{1}{(1-2^{1-r})r!}\sum_{n=1}^{\infty}\frac{1}{2^{n+1}}Y_r\left(0!H_n^{(1)}, 1!H_n^{(2)}, \ldots, (r-1)!H_n^{(r)}\right)$$

## SOME CONNECTIONS WITH THE BETA FUNCTION

Euler's beta function is defined for $\text{Re}(x) > 0$ and $\text{Re}(y) > 0$ by the integral

$$B(x, y) = \int_0^1 t^{x-1}(1-t)^{y-1}\, dt$$

and it is well known that

$$B(x, y) = \frac{\Gamma(x)\Gamma(y)}{\Gamma(x+y)}$$

Differentiating with respect to $x$ gives us

$$\frac{\partial}{\partial x}B(x, y) = \frac{\Gamma(x)\Gamma(y)}{\Gamma(x+y)}[\psi(x) - \psi(x+y)] = \int_0^1 t^{x-1}(1-t)^{y-1}\log t\, dt$$

and with $y = n$ we have

(4.1) $$\int_0^1 t^{x-1}(1-t)^{n-1}\log t\, dt = \frac{\Gamma(x)\Gamma(n)}{\Gamma(x+n)}[\psi(x) - \psi(x+n)]$$

or equivalently

(4.2) $$\int_0^1 t^{x-1}(1-t)^{n-1}\log t\, dt = -\frac{\Gamma(x)\Gamma(n)}{\Gamma(x+n)}H_n^{(1)}(x)$$

With $x = 1$ we obtain the well known integral representation for the harmonic numbers

(4.3) $$-n\int_0^1 (1-t)^{n-1}\log t\, dt = H_n^{(1)}$$

More generally we have

(4.4) $$\int_0^1 t^{x-1}(1-t)^{n-1} \log^r t \, dt = \frac{\Gamma(x)\Gamma(n)}{\Gamma(x+n)} Y_r\left(-0!H_n^{(1)}(x), 1!H_n^{(2)}(x), \ldots, (-1)^r(r-1)!H_n^{(r)}(x)\right)$$

and with $x=1$ we have

(4.5) $$\int_0^1 (1-t)^{n-1} \log^r t \, dt = \frac{1}{n} Y_r\left(-0!H_n^{(1)}, 1!H_n^{(2)}, \ldots, (-1)^r(r-1)!H_n^{(r)}\right)$$

which immediately gives us the specific value with $r=2$

(4.6) $$n\int_0^1 (1-t)^{n-1} \log^2 t \, dt = H_n^{(2)} + \left(H_n^{(1)}\right)^2$$

The following formula was also reported by Devoto and Duke [11, p.30]

$$\int_0^1 (1-t)^{n-1} \log^2 t \, dt = \frac{2}{n}\left[H_n^{(2)} + \sum_{k=1}^{n-1} \frac{H_k^{(1)}}{k+1}\right]$$

and the equivalence is readily seen by reference to Adamchik's formula [1]

(4.7) $$\sum_{k=1}^n \frac{H_k^{(1)}}{k} = \frac{1}{2}\left(H_n^{(1)}\right)^2 + \frac{1}{2} H_n^{(2)}$$

As a matter of interest, I also found formula (4.7) reported by Levenson in a 1938 volume of The American Mathematical Monthly [15] in a problem concerning the evaluation of

(4.8) $$\Gamma''(1) = \int_0^\infty e^{-x} \log^2 x \, dx = \gamma^2 + \varsigma(2)$$

We also see from (4.5) that

(4.9) $$-n\int_0^1 (1-t)^{n-1} \log^3 t \, dt = 6\left(\frac{1}{6}\left[H_n^{(1)}\right]^3 + \frac{1}{2} H_n^{(1)} H_n^{(2)} + \frac{1}{3} H_n^{(3)}\right)$$

The following formula was also reported by Devoto and Duke [11, p.30]

$$-n\int_0^1 (1-t)^{n-1} \log^3 t \, dt = 6\left[H_n^{(3)} + \sum_{k=1}^{n-1} \frac{H_k^{(1)}}{(k+1)^2} + \sum_{k=1}^{n-1} \frac{H_k^{(2)}}{k+1} + \sum_{k=1}^{n-1} \frac{1}{k+1}\sum_{j=1}^{k-2} \frac{H_j^{(1)}}{j+1}\right]$$

and its equivalence to (4.9) is readily seen by reference to Adamchik's formula [1]

$$\sum_{k=1}^{n}\frac{H_k^{(2)}}{k}+\sum_{k=1}^{n}\frac{H_k^{(1)}}{k^2}=H_n^{(3)}+H_n^{(1)}H_n^{(2)}$$

Comparing (4.5) with (3.13) we immediately see that

(4.10) $$(-1)^{r+1}n\int_0^1 (1-t)^{n-1}\log^r t\, dt = r!\sum_{k=1}^{n}\binom{n}{k}\frac{(-1)^k}{k^r}$$

which was also derived in equation (4.4.155zi) of [8].

## AN APPLICATION OF THE DIGAMMA FUNCTION

We have the classical formula for the digamma function

(5.1) $$\psi(x+a)-\psi(a)=\sum_{k=1}^{\infty}\frac{(-1)^{k+1}}{k}\frac{x(x-1)...(x-k+1)}{a(a+1)...(a+k-1)}$$

which converges for $\operatorname{Re}(x+a) > 0$. According to Raina and Ladda [19], this summation formula is due to Nörlund (see [16], [17] and also Ruben's note [20]).

Then using

$$(x)_n = x(x+1)\cdots(x+n-1) = \frac{\Gamma(x+n)}{\Gamma(x)}$$

we see that we may write

(5.2) $$\psi(a-x)-\psi(a)=-\sum_{n=1}^{\infty}\frac{1}{n}\frac{(x)_n}{(a)_n}$$

Differentiating this with respect to $x$ and using (1.22) we get

$$(-1)^r \psi^{(r)}(a-x)=-\sum_{n=1}^{\infty}\frac{1}{n}\frac{(x)_n}{(a)_n}Y_r\left(H_n^{(1)}(x),-1!H_n^{(2)}(x),...,(-1)^{r-1}(r-1)!H_n^{(r)}(x)\right)$$

Differentiating (5.2) with respect to $a$ and using

$$\frac{d}{da}\frac{1}{(a)_n}=-\frac{1}{(a)_n}H_n^{(1)}(a)$$

we obtain

$$(5.3) \quad \psi^{(r)}(a-x) - \psi^{(r)}(a) = -\sum_{n=1}^{\infty} \frac{1}{n} \frac{(x)_n}{(a)_n} Y_r\left(-H_n^{(1)}(a), 1!H_n^{(2)}(a), \ldots, (-1)^r (r-1)! H_n^{(r)}(a)\right)$$

As shown in [7] we may obtain numerous Euler sums by utilising these formulae.

From (4.4) we have

$$\frac{1}{\Gamma(n)} \int_0^1 t^{a-1}(1-t)^{n-1} \log^r t \, dt = \frac{1}{(a)_n} Y_r\left(-0! H_n^{(1)}(a), 1! H_n^{(2)}(a), \ldots, (-1)^r (r-1)! H_n^{(r)}(a)\right)$$

and substituting this in (5.3) gives us

$$\psi^{(r)}(a-x) - \psi^{(r)}(a) = -\sum_{n=1}^{\infty} \frac{(x)_n}{n!} \int_0^1 t^{a-1}(1-t)^{n-1} \log^r t \, dt$$

Differentiation with respect to $x$ results in

$$\psi^{(r+1)}(a-x) = \sum_{n=1}^{\infty} \frac{(x)_n H_n^{(1)}(x)}{n!} \int_0^1 t^{a-1}(1-t)^{n-1} \log^r t \, dt$$

and with $x = 1$ we have

$$\psi^{(r+1)}(a-1) = \sum_{n=1}^{\infty} H_n^{(1)} \int_0^1 t^{a-1}(1-t)^{n-1} \log^r t \, dt$$

$$= \int_0^1 \frac{t^{a-1} \log^r t}{1-t} \sum_{n=1}^{\infty} H_n^{(1)} (1-t)^n \, dt$$

From (2.3.1) we see that

$$\sum_{n=1}^{\infty} H_n^{(1)} (1-t)^n = -\frac{\log t}{t}$$

$$\psi^{(r+1)}(a-1) = -\int_0^1 \frac{t^{a-2} \log^{r+1} t}{1-t} \, dt$$

We may express this as

$$\psi^{(r)}(a) = -\int_0^1 \frac{t^{a-1} \log^r t}{1-t} dt$$

We note that [23, p.22]

$$\psi^{(r)}(a) = (-1)^{r+1} r! \varsigma(r+1, a)$$

which results in the well-known integral

$$(-1)^r r! \varsigma(r+1, a) = \int_0^1 \frac{t^{a-1} \log^r t}{1-t} dt$$

## SOME OTHER APPLICATIONS OF THE COMPLETE BELL POLYNOMIALS

Let us consider the function $h(x)$ with the following Maclaurin expansion

(6.1) $$\log h(x) = b_0 + \sum_{n=1}^{\infty} \frac{b_n}{n} x^n$$

and we wish to determine the coefficients $a_n$ such that

(6.2) $$h(x) = \sum_{r=0}^{\infty} a_r x^r$$

By differentiating (6.1) we obtain

$$h'(x) = h(x) \sum_{n=1}^{\infty} b_n x^{n-1} \equiv h(x) g(x)$$

From (1.13) we have

$$\frac{d^r}{dx^r} h(x) = \frac{d^r}{dx^r} e^{\log h(x)} = h(x) Y_r \left( g(x), g^{(1)}(x), \ldots, g^{(r-1)}(x) \right)$$

and in particular we have

$$\left. \frac{d^r}{dx^r} h(x) \right|_{x=0} = h(0) Y_r \left( g(0), g^{(1)}(0), \ldots, g^{(r-1)}(0) \right)$$

Using (6.2) the Maclaurin series gives us

$$a_r = \frac{1}{r!}\frac{d^r}{dx^r}h(x)\bigg|_{x=0}$$

$$a_r = \frac{1}{r!}h(0)Y_r\left(g(0), g^{(1)}(0), \ldots, g^{(r-1)}(0)\right)$$

We have

$$g^{(j)}(x) = \sum_{n=1}^{\infty} b_n (n-1)(n-2)\cdots(n-j) x^{n-1-j}$$

and thus

$$g^{(j)}(0) = j! b_{j+1}$$

Therefore we obtain

$$a_r = \frac{1}{r!} e^{b_0} Y_r\left(g(0), g^{(1)}(0), \ldots, g^{(r-1)}(0)\right)$$

$$= \frac{1}{r!} e^{b_0} Y_r\left(b_1, 1!b_2, \ldots, (r-1)!b_r\right)$$

Since $\log h(0) = \log a_0 = b_0$ we have

(6.3) $$h(x) = e^{b_0} \sum_{n=0}^{\infty} Y_n\left(b_1, 1!b_2, \ldots, (n-1)!b_n\right)\frac{x^n}{n!}$$

Then referring to (1.12)

$$\exp\left(\sum_{j=1}^{\infty} x_j \frac{t^j}{j!}\right) = \sum_{n=0}^{\infty} Y_n(x_1, \ldots, x_n)\frac{t^n}{n!}$$

we see that

$$h(x) = e^{b_0} \exp\left(\sum_{j=1}^{\infty} \frac{b_j}{j} t^j\right)$$

and this is where we started from in (6.1).

$$\log h(x) = b_0 + \sum_{n=1}^{\infty} \frac{b_n}{n} x^n$$

Multiplying (6.1) by $\alpha$ it is easily seen that

(6.4) $$h^\alpha(x) = e^{\alpha b_0} \sum_{n=0}^{\infty} Y_n(\alpha b_1, 1!\alpha b_2, \ldots, (n-1)!\alpha b_n) \frac{x^n}{n!}$$

and, in particular, with $\alpha = -1$ we obtain

(6.5) $$\frac{1}{h(x)} = e^{-b_0} \sum_{n=0}^{\infty} Y_n(-b_1, -1!b_2, \ldots, -(n-1)!b_n) \frac{x^n}{n!}$$

Differentiating (6.4) with respect to $\alpha$ would give us an expression for $h^\alpha(x) \log h(x)$.

Using (1.20) we find that

$$(r+1)! e^{-b_0} a_{r+1} = \sum_{k=0}^{r} \binom{r}{k} (r-k)! e^{-b_0} a_{r-k} k! b_{k+1}$$

giving us the recurrence relation

$$(r+1) a_{r+1} = \sum_{k=0}^{r} a_{r-k} b_{k+1}$$

or equivalently

(6.6) $$r a_r = \sum_{m=1}^{r} a_{r-m} b_m$$

Suppose that $h'(x) = h(x) g(x)$ and let $f(x) = \log h(x)$. We see that

$$f'(x) = \frac{h'(x)}{h(x)} = g(x)$$

and then using (1.9) above we have

$$\frac{d^r}{dx^r} h(x) = \frac{d^r}{dx^r} e^{\log h(x)} = h(x) Y_r\left(g(x), g^{(1)}(x), \ldots, g^{(r-1)}(x)\right)$$

Hence we have the Maclaurin expansion

$$h(x) = h(0) \sum_{r=0}^{\infty} Y_r\left(g(0), g^{(1)}(0), \ldots, g^{(r-1)}(0)\right) \frac{x^r}{r!}$$

Reference to (1.12) gives us

$$h(x) = h(0)\exp\left(\sum_{j=1}^{\infty} g^{(j-1)}(0)\frac{t^j}{j!}\right)$$

and hence we have

$$\log h(x) = \log h(0) + \sum_{j=1}^{\infty} g^{(j-1)}(0)\frac{t^j}{j!}$$

Therefore, as expected, we obtain

$$f(x) = f(0) + \sum_{j=1}^{\infty} f^{(j)}(0)\frac{t^j}{j!}$$

**Gamma function**

We note the well known series expansion

$$\log \Gamma(1+x) = -\gamma x + \sum_{n=2}^{\infty} (-1)^n \frac{\varsigma(n)}{n} x^n \quad , -1 < x \leq 1$$

and hence we have

(6.7) $$\Gamma(1+x) = \sum_{n=0}^{\infty} Y_n\left(-0!\varsigma(1), 1!\varsigma(2), \ldots, (-1)^n (n-1)!\varsigma(n)\right)\frac{x^n}{n!}$$

and

(6.8) $$\frac{1}{\Gamma(1+x)} = \sum_{n=0}^{\infty} Y_n\left(0!\varsigma(1), -1!\varsigma(2), \ldots, (-1)^{n+1}(n-1)!\varsigma(n)\right)\frac{x^n}{n!}$$

where $\varsigma(1)$ is defined as equal to $\gamma$.

Differentiating (6.7) we get

$$\Gamma^{(m)}(1+x) = \sum_{n=1}^{\infty} Y_n\left(-0!\varsigma(1), 1!\varsigma(2), \ldots, (-1)^n(n-1)!\varsigma(n)\right)\frac{n(n-1)\ldots(n-m+1)x^{n-m}}{n!}$$

and hence, letting $x = 0$, we have the $m$ th derivative of the gamma function

(6.9) $$\Gamma^{(m)}(1) = Y_n\left(-0!\varsigma(1), 1!\varsigma(2), \ldots, (-1)^m(m-1)!\varsigma(m)\right)$$

where we again designate $\varsigma(1) = \gamma$.

It is reported as an exercise in Apostol's book [2, p.303] that $\Gamma^{(n)}(1)$ has the same sign as $(-1)^n$ and a derivation is contained in [9].

An alternative proof is shown below. Since $\Gamma'(x) = \Gamma(x)\psi(x)$ we have

$$\Gamma^{(m)}(x) = \Gamma(x)Y_m\left(\psi(x), \psi^{(1)}(x), ..., \psi^{(m-1)}(x)\right)$$

and using (3.7)

$$\psi^{(p)}(x) = (-1)^{p+1} p!\varsigma(p+1, x)$$

we may express $\Gamma^{(m)}(x)$ in terms of $\psi(x)$ and the Hurwitz zeta functions.

$$\Gamma^{(m)}(x) = \Gamma(x)Y_m\left(\psi(x), 1!\varsigma(2, x), ..., (-1)^m (m-1)!\varsigma(m, x)\right)$$

From the definition of the (exponential) complete Bell polynomials we have

$$Y_m(ax_1, a^2 x_2, ..., a^m x_m) = a^m Y_r(x_1, ..., x_m)$$

and thus with $a = -1$ we have

$$Y_m(-x_1, x_2, ..., (-1)^m x_m) = (-1)^m Y_m(x_1, ..., x_m)$$

Hence we see that

$$\Gamma^{(m)}(x) = (-1)^m \Gamma(x)Y_m\left(-\psi(x), 1!\varsigma(2, x), ..., (m-1)!\varsigma(m, x)\right)$$

It is well known that $\psi(x)$ is negative in the interval $[0, \alpha)$ where $\alpha > 0$ is the unique solution of $\psi(\alpha) = 0$ and hence we see that $-\psi(x)$ is positive in that interval. Since $Y_m(x_1, ..., x_m) > 0$ when all of the arguments are positive, we deduce that $\Gamma^{(m)}(x)$ has the same sign as $(-1)^m$ when $x \in [0, \alpha)$.

We also note that

$$Y_m(x_1, -x_2, ..., (-1)^{m+1} x_m) = (-1)^m Y_m(-x_1, ..., -x_m)$$

but no discernable sign pattern emerges here. We may note that

$$Y_n\left(b_1, -1!b_2, \ldots, (-1)^{n+1}(n-1)!b_n\right) = (-1)^n Y_n\left(-b_1, -1!b_2, \ldots, -(n-1)!b_n\right)$$

$$= (-1)^n e^{-b_0} \frac{d^n}{dx^n} \frac{1}{h(x)}\bigg|_{x=0}$$

We also have the expansion [23, p.159]

$$\log \Gamma(a+x) = \log \Gamma(a) + \psi(a)x + \sum_{n=2}^{\infty} \frac{(-1)^n \varsigma(n,a)}{n} x^n$$

and hence we deduce that

$$\Gamma(a+x) = \Gamma(a) \sum_{n=0}^{\infty} Y_n\left(0!\psi(a), 1!\varsigma(2,a), \ldots, (-1)^n (n-1)!\varsigma(n,a)\right) \frac{x^n}{n!}$$

**Barnes double gamma function**

We note the Barnes double gamma function $G(x)$ defined by [23, p.25]

$$G(1+x) = (2\pi)^{x/2} \exp\left[-\frac{1}{2}(\gamma x^2 + x^2 + x)\right] \prod_{k=1}^{\infty} \left\{\left(1 + \frac{x}{k}\right)^k \exp\left(\frac{x^2}{2k} - x\right)\right\}$$

The following identity was originally derived by Srivastava [23, p.210] in 1988

$$\log G(1+x) = \frac{1}{2}[\log(2\pi) - 1]x - \frac{1}{2}(1+\gamma)x^2 + \sum_{k=3}^{\infty} (-1)^{k+1} \frac{\varsigma(k-1)}{k} x^k$$

and from this we may determines a series expansion for $1/G(x+1)$ in terms involving the exponential Bell polynomials. Using () we have

$$\frac{1}{G(1+x)} = \sum_{n=0}^{\infty} Y_n\left(-c_1, 1!c_2, \ldots, (-1)^n (n-1)!c_n\right) \frac{x^n}{n!}$$

where

$$c_1 = \frac{1}{2}[\log(2\pi) - 1] \qquad c_2 = (1+\gamma) \qquad c_n = \varsigma(n-1) \text{ for } n \geq 3$$

and we note that $c_n > 0$ for all $n \geq 1$.

We then see that

$$\left.\frac{d^m}{dx^m}\frac{1}{G(1+x)}\right|_{x=0} = Y_m\left(-c_1, 1!c_2, \ldots, (-1)^m(m-1)!c_m\right)$$

and noting (1.10.2) we determine that

$$\left.\frac{d^m}{dx^m}\frac{1}{G(1+x)}\right|_{x=0} = (-1)^m Y_m\left(c_1, 1!c_2, \ldots, (m-1)!c_m\right)$$

so that the derivatives have the same sign as $(-1)^m$. In fact $\frac{d^m}{dx^m}\frac{1}{G(1+x)}$ has the same sign as $(-1)^m$ for all $x \geq 0$.

For example we have

$$\left.\frac{d}{dx}\frac{1}{G(1+x)}\right|_{x=0} = -\left.\frac{G'(1+x)}{G^2(1+x)}\right|_{x=0} = -G'(1)$$

and hence we have

$$G'(1) = \frac{1}{2}\left[\log(2\pi) - 1\right]$$

**Bernoulli numbers**

The Bernoulli numbers are defined by the generating function

$$\frac{t}{e^t - 1} = \sum_{n=0}^{\infty} \frac{B_n}{n!} t^n$$

and thus

$$\frac{1}{e^t - 1} = \sum_{n=0}^{\infty} \frac{B_n}{n!} t^{n-1}$$

Integration gives us

$$\int_a^x \frac{dt}{e^t - 1} = \int_a^x \frac{e^{-t} dt}{1 - e^{-t}} = \log(1 - e^{-x}) - \log(1 - e^{-a})$$

$$= \log x - \log a + \sum_{n=1}^{\infty} \frac{B_n}{n\, n!}[x^n - a^n]$$

and hence we obtain

$$\log(1-e^{-x}) - \log(1-e^{-a}) = \log x - \log a + \sum_{n=1}^{\infty} \frac{B_n}{n\,n!}[x^n - a^n]$$

We may express this as

$$\log \frac{(1-e^{-x})}{x} = \log \frac{(1-e^{-a})}{a} + \sum_{n=1}^{\infty} \frac{B_n}{n\,n!}[x^n - a^n]$$

Using L'Hôpital's rule we see that

$$\lim_{a \to 0} \log \frac{(1-e^{-a})}{a} = 0$$

and with $x \to -x$ we obtain

$$\log\left(\frac{x}{e^x - 1}\right) = \sum_{n=1}^{\infty} \frac{(-1)^{n+1} B_n}{n\,n!} x^n$$

which was proved by Ramanujan for $|x| < 2\pi$ [3, p.119].

Hence we have

$$\frac{x}{e^x - 1} = \sum_{n=0}^{\infty} Y_n\left(\frac{B_1}{1}, -\frac{B_2}{2}, \ldots, (-1)^{n+1}\frac{B_n}{n}\right) \frac{x^n}{n!}$$

and we therefore see that

$$\sum_{n=0}^{\infty} \frac{B_n}{n!} x^n = \sum_{n=0}^{\infty} Y_n\left(\frac{B_1}{1}, -\frac{B_2}{2}, \ldots, (-1)^{n+1}\frac{B_n}{n}\right) \frac{x^n}{n!}$$

Equating coefficients gives us the recurrence relation

(6.10) $$B_n = Y_n\left(\frac{B_1}{1}, -\frac{B_2}{2}, \ldots, (-1)^{n+1}\frac{B_n}{n}\right)$$

Employing (1.20) we obtain

(6.11) $$B_{r+1} = \sum_{k=0}^{r} (-1)^k \binom{r}{k} \frac{B_{r-k} B_{k+1}}{k+1}$$

**Values of the Riemann zeta function $\varsigma(2n)$**

This section is based on a paper by Snowden [22]. We have the well-known infinite product

$$\sin \pi x = \pi x \prod_{k=1}^{\infty}\left(1-\frac{x^2}{k^2}\right)$$

and we have for $|x|<1$

$$\log\left(\frac{\sin \pi x}{\pi x}\right) = \sum_{k=1}^{\infty}\log\left(1-\frac{x^2}{k^2}\right) = -\sum_{n=1}^{\infty}\frac{\varsigma(2n)}{n}x^{2n}$$

$$= -2\sum_{n=1}^{\infty}\frac{b_n}{n}x^n$$

where $b_{2n} = 2\varsigma(2n)$, $b_{2n+1} = 0$.

Therefore we have

$$\frac{\sin \pi x}{\pi x} = \sum_{n=0}^{\infty} Y_n\left(-b_1, -1!b_2, \ldots, -(n-1)!b_n\right)\frac{x^n}{n!}$$

We also have the Maclaurin series

$$\frac{\sin \pi x}{\pi x} = \sum_{n=0}^{\infty}(-1)^n \pi^{2n}\frac{x^{2n}}{(2n+1)!}$$

and equating coefficients gives us

$$\frac{(-1)^n (2n)! \pi^{2n}}{(2n+1)!} = Y_{2n}\left(-b_1, -1!b_2, \ldots, -(2n-1)!b_{2n}\right)$$

$$0 = Y_{2n+1}\left(-b_1, -1!b_2, \ldots, -(2n)!b_{2n+1}\right)$$

For example, with $n=1$ we obtain Euler's formula

$$-\frac{\pi^2}{3} = Y_2\left(0, -2\varsigma(2)\right) \text{ or } \varsigma(2) = \frac{\pi^2}{6}$$

and with $n=2$ we have

$$\frac{\pi^4}{60} = \varsigma^2(2) - \varsigma(4)$$

We obtain from (1.20)

$$Y_{2n} = \sum_{k=0}^{2n-1} \binom{2n-1}{k} Y_{2n-1-k}\, x_{k+1}$$

$$= \sum_{k=0}^{n} \binom{2n-1}{2k} Y_{2n-1-2k}\, x_{2k+1} + \sum_{k=1}^{n} \binom{2n-1}{2k-1} Y_{2n-2k}\, x_{2k}$$

and, because $b_{2n+1} = 0$, we have in our specific case

$$Y_{2n} = \sum_{k=1}^{n} \binom{2n-1}{2k-1} Y_{2n-2k}\, x_{2k}$$

This results in

$$\frac{(-1)^n (2n)!\, \pi^{2n}}{(2n+1)!} = -2 \sum_{k=1}^{n} \binom{2n-1}{2k-1} \frac{(-1)^{n-k}(2n-2k)!\, \pi^{2n-2k}(2k-1)!\, \varsigma(2k)}{(2n-2k+1)!}$$

which simplifies to

$$\frac{(2n)!}{(2n-1)!(2n+1)!} = 2 \sum_{k=1}^{n} \frac{(-1)^{k+1} \pi^{-2k} \varsigma(2k)}{(2n-2k+1)!}$$

Then substituting the well known relation [23, p.98]

$$\varsigma(2k) = \frac{(-1)^{k+1} 2^{2k} \pi^{2k} B_{2k}}{2(2k)!}$$

we obtain

$$\frac{(2n)!}{(2n-1)!(2n+1)!} = \sum_{k=1}^{n} \frac{2^{2k} B_{2k}}{(2n-2k+1)!(2k)!}$$

or equivalently

$$\frac{[(2n)!]^2}{(2n-1)!(2n+1)!} = \sum_{k=1}^{n} \binom{2n}{2k} \frac{2^{2k} B_{2k}}{(2n-2k+1)}$$

**REFERENCES**

[1]  V.S.Adamchik, On Stirling Numbers and Euler Sums.

    J. Comput. Appl. Math.79, 119-130, 1997.
    http://www-2.cs.cmu.edu/~adamchik/articles/stirling.htm

[2] T.M. Apostol, Mathematical Analysis, Second Ed., Addison-Wesley Publishing Company, Menlo Park (California), London and Don Mills (Ontario), 1974.

[3] B.C. Berndt, Ramanujan's Notebooks. Part I, Springer-Verlag, 1985.

[4] M.W. Coffey, A set of identities for a class of alternating binomial sums arising in computing applications. 2006.
arXiv:math-ph/0608049v1

[5] C.A. Charalambides, Enumerative Combinatorics.
Chapman & Hall/CRC, 2002.

[6] L. Comtet, Advanced Combinatorics, Reidel, Dordrecht, 1974.

[7] D.F. Connon, Some series and integrals involving the Riemann zeta function, binomial coefficients and the harmonic numbers. Volume I, 2007.
arXiv:0710.4022 [pdf]

[8] D.F. Connon, Some series and integrals involving the Riemann zeta function, binomial coefficients and the harmonic numbers. Volume II(a), 2007.
arXiv:0710.4023 [pdf]

[9] D.F. Connon, Some integrals involving the Stieltjes constants. 2009.
arXiv:0902.2188 [pdf]

[10] M.A. Coppo, La formule d'Hermite revisitée. 2003.
http://math1.unice.fr/~coppo/

[11] A. Devoto and D.W. Duke, Table of integrals and formulae for Feynman diagram calculations. Florida State University, FSU-HEP-831003, 1983.
http://www.csit.fsu.edu/~dduke/integrals.htm

[12] P. Flajolet and R. Sedgewick, Mellin Transforms and Asymptotics: Finite Differences and Rice's Integrals. Theor. Comput. Sci. 144, 101-124, 1995.
Mellin transforms and asymptotics: Finite differences and Rice's integrals
psu.edu [pdf]

[13] R.L. Graham, D.E. Knuth and O. Patashnik, Concrete Mathematics. Second Ed. Addison-Wesley Publishing Company, Reading, Massachusetts, 1994.

[14] K.S. Kölbig, The complete Bell polynomials for certain arguments in terms of Stirling numbers of the first kind.
J. Comput. Appl. Math. 51 (1994) 113-116.
Also available electronically at:

A relation between the Bell polynomials at certain arguments and a Pochhammer symbol.
CERN/Computing and Networks Division, CN/93/2, 1993 http://cdsweb.cern.ch/

[15] M.E. Levenson, J.F. Locke and H. Tate, Amer. Math. Monthly, 45, 56-58, 1938.

[16] N.E. Nörlund, Vorlesungen über Differenzenrechnung. Chelsea, 1954.
http://dz-srv1.sub.uni-goettingen.de/cache/browse/AuthorMathematicaMonograph,WorkContainedN1.html

[17] N.E. Nörlund, Leçons sur les séries d'interpolation.
Paris, Gauthier-Villars, 1926.

[18] G. Póyla and G. Szegö, Problems and Theorems in Analysis, Vol. I.
Springer-Verlag, New York 1972.

[19] R.K. Raina and R.K. Ladda, A new family of functional series relations involving digamma functions.
Ann. Math. Blaise Pascal, Vol. 3, No. 2, 1996, 189-198.
http://www.numdam.org/item?id=AMBP_1996__3_2_189_0

[20] H. Ruben, A Note on the Trigamma Function.
Amer. Math. Monthly, 83, 622-623, 1976.

[21] L.-C. Shen, Remarks on some integrals and series involving the Stirling numbers and $\varsigma(n)$. Trans. Amer. Math. Soc. 347, 1391-1399, 1995.

[22] A. Snowden, Collection of Mathematical Articles. 2003.
http://www-users.math.umd.edu/~asnowden/math-cont/dorfman.pdf

[23] H.M. Srivastava and J. Choi, Series Associated with the Zeta and Related Functions. Kluwer Academic Publishers, Dordrecht, the Netherlands, 2001.

[24] H.S. Wilf, The asymptotic behaviour of the Stirling numbers of the first kind.
Journal of Combinatorial Theory, Series A, 64, 344-349, 1993.
http://citeseer.ist.psu.edu/520553.html

Donal F. Connon
Elmhurst
Dundle Road
Matfield, Kent TN12 7HD
dconnon@btopenworld.com